\documentclass[12pt]{article}

\usepackage{graphicx}

\def\be{\begin{equation}}
\def\ee{\end{equation}}
\newcommand{\kk}[2]{\frac{#1}{#2}}

\newcommand{\ff}[1]{{\bf  #1}}
\def\a{\alpha}
\def\b{\beta}

\def\e{\epsilon}
\def\lam{\lambda}
\def\s{\qquad}
\def\ra{\rightarrow}
\def\={\approx}
\def\8{{\infty}}
\def\x{\ff{x}}
\def\vcode#1#2#3#4{\begin{figure}
\begin{center}
\begin{minipage}[c]{#1\textwidth}
{{\small #2 \hrule \vspace{5pt} {\bf begin}  \\   %
{\it #3} \\ {\bf end} \vspace{5pt} \hrule }}
\end{minipage}
\caption{#4}
\end{center}   \end{figure}} %%
\begin{document}

\title{Firefly Algorithm, L\'evy Flights and Global Optimization}

\author{Xin-She Yang \\
Department of Engineering, University of Cambridge, \\
Trumpington Street, Cambridge CB2 1PZ, UK }

\date{}

\maketitle

%% Abstract

\abstract{
Nature-inspired algorithms such as Particle Swarm Optimization and Firefly Algorithm
are among the most powerful algorithms for optimization.  In this paper, we intend
to formulate a new metaheuristic algorithm by combining L\'evy flights with
the search strategy via the Firefly Algorithm. Numerical studies and results
suggest that the proposed L\'evy-flight firefly algorithm is
superior to existing metaheuristic algorithms. Finally
implications for further research and wider applications will be discussed. \\
}

\noindent {\bf Citation detail:} X.-S. Yang,``Firefly algorithm, L\'evy flights and global optimization",
in: {\it Research and Development in Intelligent Systems XXVI} (Eds M. Bramer, R. Ellis, M. Petridis),
Springer London, pp. 209-218 (2010). 

%% Begin of Main Text %%

\section{Introduction}

Nature-inspired metaheuristic algorithms are becoming powerful in solving modern
global optimization problems \cite{Baeck,Bod,Deb,Gold,Ken2,Yang,Yang2},
especially for the NP-hard optimization such
as the travelling salesman problem. For example,
particle swarm optimization (PSO) was developed by Kennedy and
Eberhart in 1995 \cite{Ken,Ken2}, based on the swarm behaviour such
as fish and bird schooling in nature. It has now been applied
to find solutions for many optimization applications.
Another example is the Firefly Algorithm developed by the author
\cite{Yang} which has demonstrated promising superiority over
many other algorithms. The search strategies in these multi-agent
algorithms are controlled randomization, efficient local search and selection of the
best solutions. However, the randomization typically uses
uniform distribution or Gaussian distribution.

On the other hand, various studies have shown that
flight behaviour of many animals and insects has demonstrated
the typical characteristics of L\'evy flights
\cite{Brown,Reynolds,Pav,Pav2}. A recent study by Reynolds and Frye shows that
fruit flies or {\it Drosophila melanogaster}, explore their landscape using a series
of straight flight paths punctuated by a sudden $90^{\o}$ turn, leading to
a L\'evy-flight-style intermittent scale free search pattern.
Studies on human behaviour such as the Ju/'hoansi hunter-gatherer foraging
patterns also show the typical feature of L\'evy flights.
Even light can be related to L\'evy flights \cite{Barth}.
Subsequently, such behaviour has been applied to
optimization and optimal search, and preliminary results show its
promising capability \cite{Pav,Reynolds,Shles,Shles2}.

This paper aims to formulate a new L\'evy-flight Firefly Algorithm (LFA)
and to provide the comparison study of the LFA with PSO and
other relevant algorithms. We will first outline
the firefly algorithms, then formulate the L\'evy-flight FA
and finally give the comparison about the performance of these algorithms.
The LFA optimization seems more promising than particle swarm optimization in the sense
that LFA converges more quickly and deals with global optimization more naturally.
In addition, particle swarm optimization is just a special class of
the LFA as we will demonstrate this in this paper.

\section{Firefly Algorithm}

\subsection{Behaviour of Fireflies}

The flashing light of fireflies is an amazing sight in the summer sky in
the tropical and temperate regions. There are about two thousand
firefly species, and most fireflies produce short and rhythmic flashes.
The pattern of flashes is often unique for a particular species.
The flashing light is produced by a process of bioluminescence, and
the true functions of such signaling systems are still debating. However,
two fundamental functions of such flashes are to attract mating partners
(communication), and to attract potential prey. In addition,
flashing may also serve as a protective warning mechanism.
The rhythmic flash, the rate of flashing and the amount of time
form part of the signal system that brings both sexes together.
Females respond to a male's unique pattern of flashing
in the same species, while in some species such as {\it photuris},
female fireflies  can mimic the mating
flashing pattern of other species so as to lure and eat the male
fireflies who may mistake the flashes as a potential suitable mate.

The flashing light can be formulated in such a way that it is associated with
the objective function to be optimized, which makes it possible to formulate
new optimization algorithms. In the rest of this paper, we will first outline
the basic formulation of the Firefly Algorithm (FA) and then discuss
the implementation as well as analysis in detail.

\subsection{Firefly Algorithm}

Now we can idealize some of the flashing characteristics of fireflies
so as to develop firefly-inspired algorithms. For simplicity in
describing our Firefly Algorithm (FA), we now use the following
three idealized rules: 1) all fireflies are unisex so that one firefly will be attracted to other fireflies
regardless of their sex; 2) Attractiveness is proportional to their brightness,
thus for any two flashing fireflies, the less brighter one will move towards the brighter
one. The attractiveness is proportional to the brightness and they both
decrease as their distance increases. If there is
no brighter one than a particular firefly, it will move randomly;
3) The brightness of a firefly is affected or  determined by the landscape of the
objective function. For a  maximization problem, the brightness can simply be proportional
to the value of the objective function. Other forms of brightness can be defined in a similar
way to the fitness function in genetic algorithms or the bacterial
foraging algorithm (BFA) \cite{Gazi,Passino}.

In the firefly algorithm, there are two important issues: the variation of
light intensity and formulation of the attractiveness.
For simplicity, we can always assume that the attractiveness
of a firefly is determined by its brightness which in turn is associated with
the encoded objective function.

In the simplest case for maximum optimization problems,
the brightness $I$ of a firefly at a particular location $\x$ can be chosen
as $I(\x) \propto f(\x)$. However, the attractiveness $\b$ is relative, it should be
seen in the eyes of the beholder or judged by the other fireflies. Thus, it will
vary with the distance $r_{ij}$ between firefly $i$ and firefly $j$. In addition,
light intensity decreases with the distance from its source, and light is also
absorbed in the media,
so we should allow the attractiveness to vary with the degree of absorption.
In the simplest form, the light intensity $I(r)$ varies according to the inverse
square law $I(r)=\kk{I_s}{r^2}$ where $I_s$ is the intensity at the source.
For a given medium with a fixed light absorption coefficient $\gamma$,
the light intensity $I$ varies with the distance $r$. That is
\be I=I_0 e^{-\gamma r}, \ee
where $I_0$ is the original light intensity.

As a firefly's attractiveness is proportional to
the light intensity seen by adjacent fireflies, we can now define
the attractiveness $\b$ of a firefly by
\be \b = \b_0 e^{-\gamma r^2}, \label{att-equ-100} \ee
where $\b_0$ is the attractiveness at $r=0$.

\section{L\'evy-Flight Firefly Algorithm}

If we combine the three idealized rules with the characteristics of
L\'evy flights, we can formulate a new L\'evy-flight Firefly Algorithm (LFA)
which can be summarized as the pseudo code shown in Fig. \ref{fa-fig-100}.

\vcode{0.9}{{\sf L\'evy-Flight Firefly Algorithm}} {
\indent \quad Objective function $f(\x), \s \x=(x_1, ..., x_d)^T$ \\
\indent \quad Generate initial population of fireflies $\x_i \; (i=1,2,...,n)$ \\
\indent \quad Light intensity $I_i$ at $\x_i$ is determined by $f(\x_i)$  \\
\indent \quad Define light absorption coefficient $\gamma$ \\
\indent \quad {\bf while} ($t<$MaxGeneration) \\
\indent \quad {\bf for} $i=1:n$ all $n$ fireflies \\
\indent \qquad {\bf for} $j=1:i$ all $n$ fireflies \\
\indent \qquad \qquad {\bf if} ($I_j>I_i$) \\
\indent \qquad \qquad Move firefly $i$ towards $j$ in d-dimension via L\'evy flights  \\
\indent \qquad \qquad {\bf end if} \\
\indent \qquad \qquad Attractiveness varies with distance $r$ via $\exp[-\gamma r]$ \\
\indent \qquad \qquad Evaluate new solutions and update light intensity \\
\indent \qquad {\bf end for }$j$ \\
\indent \quad {\bf end for }$i$ \\
\indent \quad Rank the fireflies and find the current best  \\
\indent \quad {\bf end while} \\
\indent \quad Postprocess results and visualization }{Pseudo code of the L\'evy-Flight Firefly Algorithm (LFA).
\label{fa-fig-100} }

In the implementation, the actual form of attractiveness function $\b(r)$ can be
any monotonically decreasing functions such as the following generalized form
\be \b(r) =\b_0 e^{-\gamma r^m}, \s (m \ge 1). \ee
For a fixed $\gamma$, the characteristic length becomes
$\Gamma=\gamma^{-1/m} \ra 1$ as $m \ra \infty$.  Conversely,
for a given length scale $\Gamma$ in an optimization problem, the parameter $\gamma$
can be used as a typical initial value. That is $\gamma =\kk{1}{\Gamma^m}$.

The distance between any two fireflies $i$ and $j$ at $\x_i$ and $\x_j$, respectively, is
the Cartesian distance
\be r_{ij}=||\x_i-\x_j|| =\sqrt{\sum_{k=1}^d (x_{i,k} - x_{j,k})^2}, \ee
where $x_{i,k}$ is the $k$th component of the spatial coordinate $\x_i$ of $i$th
firefly. For other applications such as scheduling, the distance can be time delay or
any suitable forms.

The movement of a firefly $i$ is attracted to another more attractive (brighter)
firefly $j$ is determined by
\be \x_i =\x_i + \b_0 e^{-\gamma r^2_{ij}} (\x_j-\x_i) + \a \; \textrm{sign}[{\rm rand}-\kk{1}{2}]
\oplus \textrm{L\'evy}, \ee
where the second term is due to the attraction while the third term
is randomization via L\'evy flights with $\a$ being the randomization parameter.
The product $\oplus$ means entrywise multiplications.
The {\sf sign[rand-$\kk{1}{2}$]} where {\sf rand} $\in [0,1]$
essentially provides a random sign or direction while the random step length
is drawn from a L\'evy distribution
\be \textrm{L\'evy} \sim u = t^{-\lam}, \s (1 < \lam \le 3), \ee
which has an infinite variance with an infinite mean.
Here the steps of firefly motion is essentially
a random walk process with a power-law step-length distribution with a heavy tail.

\subsection{Choice of Parameters}

For most cases in our implementation,
we can take $\b_0=1$, $\a \in [0,1]$, $\gamma=1$, and $\lam=1.5$.
In addition, if the scales vary significantly in different dimensions such as $-10^5$ to $10^5$ in
one dimension while, say, $-0.001$ to $0.01$ along the other, it is a good idea to
replace $\a$ by $\a S_k$ where the scaling parameters $S_k (k=1,...,d)$ in
the $d$ dimensions should be determined by the actual scales of the problem of interest.

The parameter $\gamma$ now characterizes the variation of the attractiveness,
and its value is crucially important in determining the speed of the convergence
and how the FA algorithm behaves. In theory, $\gamma \in [0,\infty)$, but in
practice, $\gamma=O(1)$ is determined by the characteristic length $\Gamma$ of the
system to be optimized. Thus, in most applications, it typically varies
from $0.01$ to $100$.

\subsection{Asymptotic Cases}

There are two important limiting cases when $\gamma \ra 0$ and $\gamma \ra \infty$.
For $\gamma \ra 0$, the attractiveness is constant $\b=\b_0$ and $\Gamma \ra \infty$,
this is equivalent to say that the light intensity does not decrease in an idealized sky.
Thus, a flashing firefly can be seen anywhere in the domain. Thus, a single (usually global)
optimum can easily be reached. This corresponds to a special case of particle
swarm optimization (PSO) discussed earlier. Subsequently, the efficiency of this special case
is the same as that of PSO.

On the other hand, the limiting case $\gamma \ra \infty$ leads to
$\Gamma \ra 0$ and  $\b(r) \ra \delta(r)$ (the Dirac delta function), which means that
the attractiveness is almost zero in the sight of other fireflies or the fireflies are short-sighted. This is equivalent
to the case where the fireflies fly in a very foggy region randomly. No other fireflies can be
seen, and each firefly roams in a completely random way. Therefore, this corresponds
to the completely random search method.

As the L\'evy-flight firefly algorithm is usually in somewhere between these two extremes, it is possible
to adjust the parameters $\gamma$, $\lambda$ and $\alpha$ so that it can outperform both the random search
and PSO. In fact, LFA can find the global optima as well as all the local optima
simultaneously in a very effective manner.

\section{Simulations and Results}

\subsection{Validation}

In order to validate the proposed algorithm, we have implemented it in
Matlab.  In our simulations,  the values
of the parameters are $\a=0.2$, $\gamma=1$, $\lam=1.5$, and $\b_0=1$.
\begin{figure}
 \centerline{\includegraphics[width=4in,height=3in]{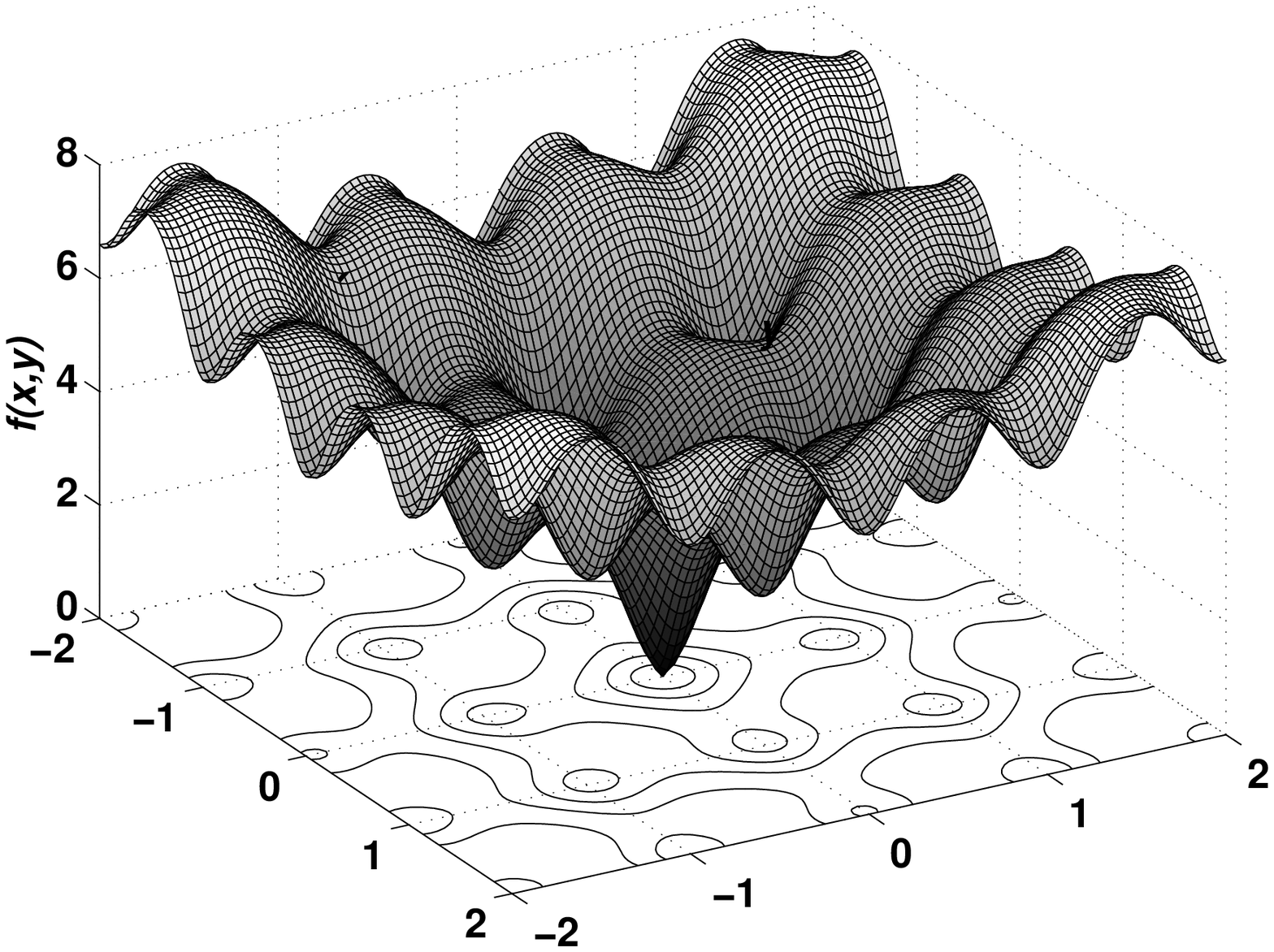} }
\caption{Ackley function for two independent variables
with a global minimum $f_*=0$ at $(0,0)$. \label{Yang-fig1} }
\end{figure}
As an example, we now use the LFA to find the
global optimum of the Ackley function
\be f(\x)=-20 \exp\Big[-\kk{1}{5} \sqrt{\kk{1}{d} \sum_{i=1}^d x_i^2}\Big] - \exp[\kk{1}{d} \sum_{i=1}^d \cos (2 \pi x_i)]
+20 +e, \ee
which has a global minimum $f_*=0$ at $(0,0,...,0)$. The 2D Ackley function is
shown in Fig. \ref{Yang-fig1}, and this global minimum can be found
after about 200 evaluations
for 40 fireflies after 5 iterations as shown in Fig. \ref{Yang-fig2}.

Now let us use the LFA to find the optima of some tougher
test functions. For example, the author introduced a forest function \cite{Yang3}
\be f(\x)=\Big( \sum_{i=1}^d |x_i| \Big) \exp\Big[- \sum_{i=1}^d \sin (x_i^2) \Big],
\;\;\; -2 \pi \le x_i \le 2 \pi, \ee
which has a global minimum $f_*=0$ at $(0,0,...,0)$. The 2D Yang's forest function
is shown in Fig. \ref{Yang-fig3}. However, an important feature of this
test function is that it is non-smooth
and its derivative is not well defined at the optima $(0,0,...,0)$
as shown in Fig. \ref{Yang-fig4}.

\subsection{Comparison of LFA with PSO and GA}

\begin{figure}
 \centerline{\includegraphics[height=3in,width=3in]{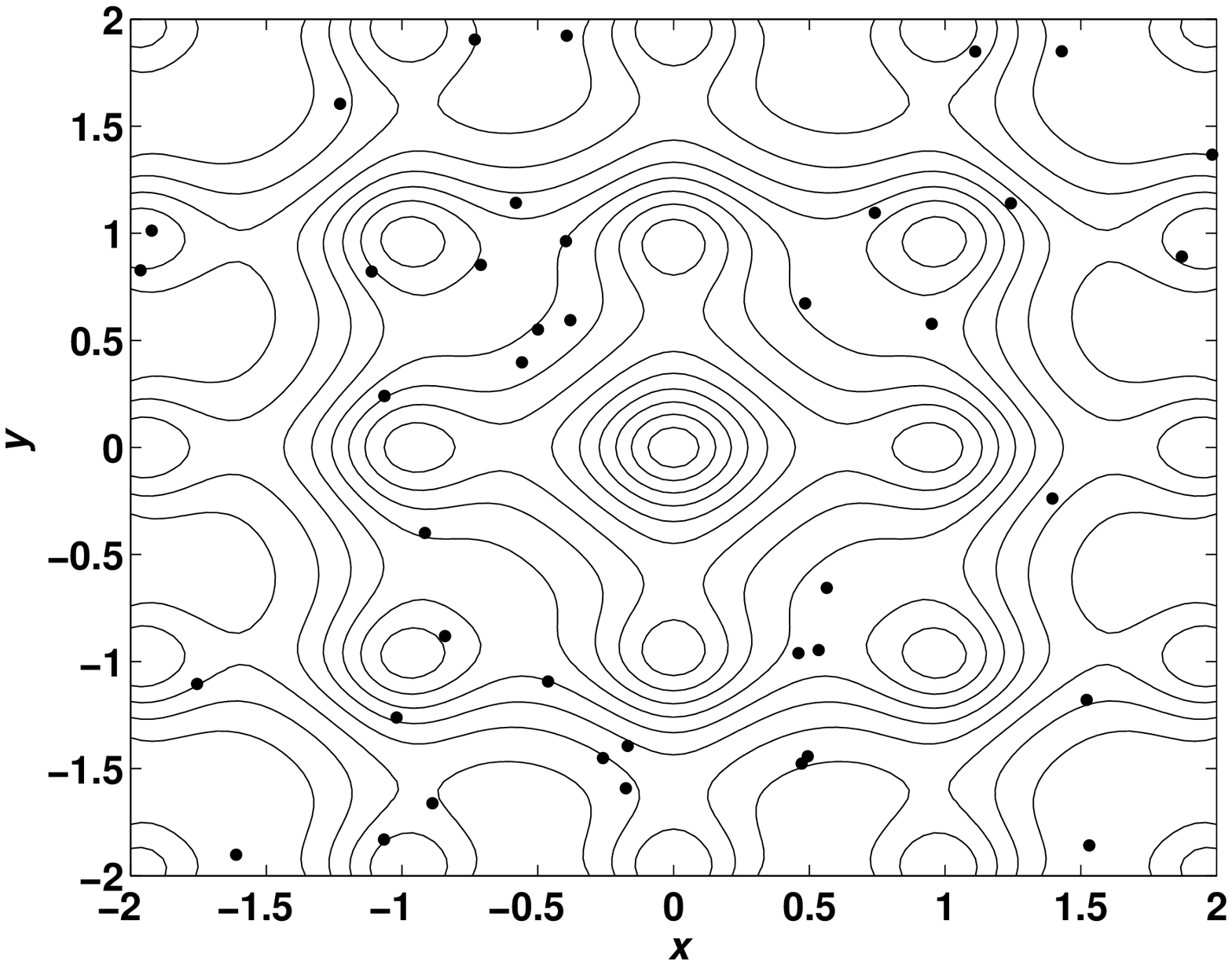}
 \includegraphics[height=3in,width=3in]{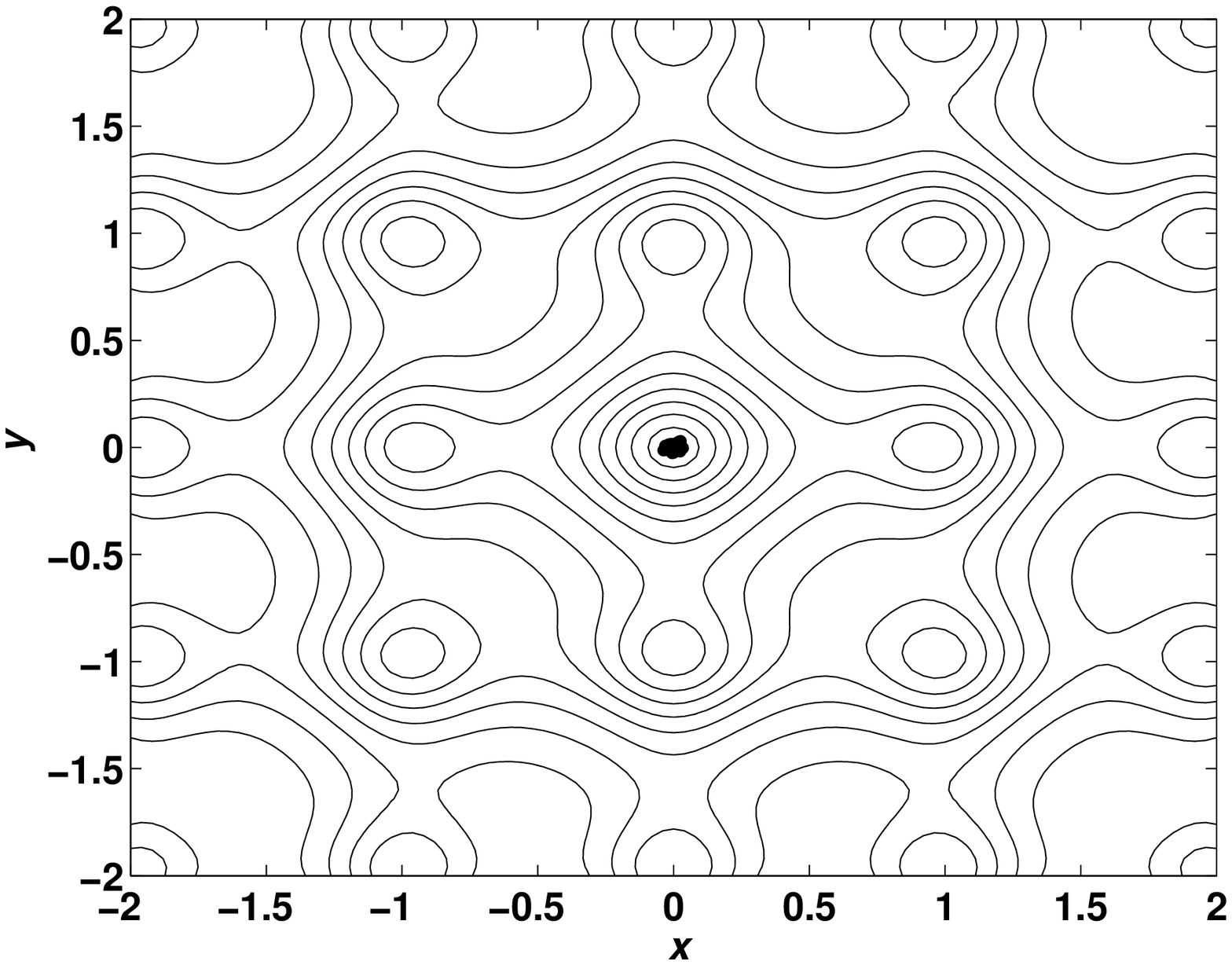}}
\vspace{-5pt}
\caption{The initial locations of the 40 fireflies (left) and their
locations after $5$ iterations (right). \label{Yang-fig2} }
\end{figure}

Various studies show that PSO algorithms can
outperform genetic algorithms (GA) \cite{Gold} and other conventional algorithms for
solving many optimization problems. This is
partially due to that fact that the broadcasting ability of the current
best estimates gives better and quicker convergence towards the
optimality. A general framework for evaluating statistical performance of evolutionary algorithms
has been discussed in detail by Shilane et al. \cite{Shilane}.  Various test functions
for optimization algorithms
have been developed over many years, and  a relatively comprehensive review of these
test functions can be found in \cite{Baeck}.

Now we will compare the LFA with PSO, and genetic
algorithms for various standard test functions. We will use the same
population size of $n=40$ for all algorithms in all our simulations.
The PSO used is the standard version without any inertia function, while
the implemented genetic algorithm  has a mutation probability of $0.05$ and a crossover
probability of $0.95$ without use of elitism.
After implementing these algorithms using
Matlab, we have carried out extensive simulations and each algorithm has been
run at least 100 times so as to carry out meaningful statistical analysis.
The algorithms stop when the variations of function values are less than
a given tolerance $\e \le 10^{-5}$. The results are
summarized in the following table (see Table 1) where the global optima
are reached. The numbers are in the
format: average number of evaluations (success rate), so
$6922 \pm 537 (98\%)$ means that the average number (mean) of function
evaluations is 6922 with a standard deviation of 537. The success rate
of finding the global optima for this algorithm is $98\%$.

\begin{figure}
 \centerline{\includegraphics[width=3in,height=3in]{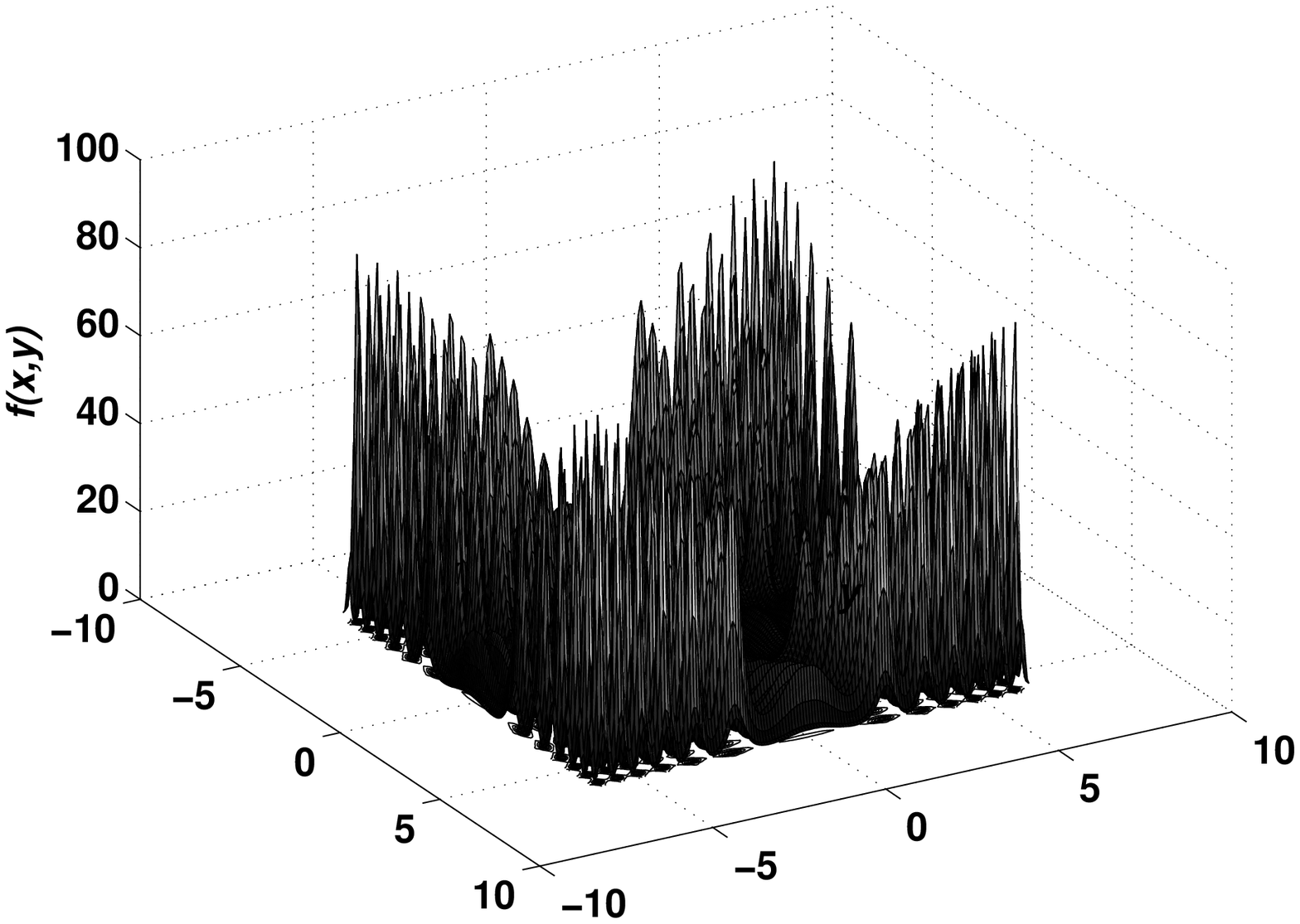}
 \includegraphics[width=3in,height=3in]{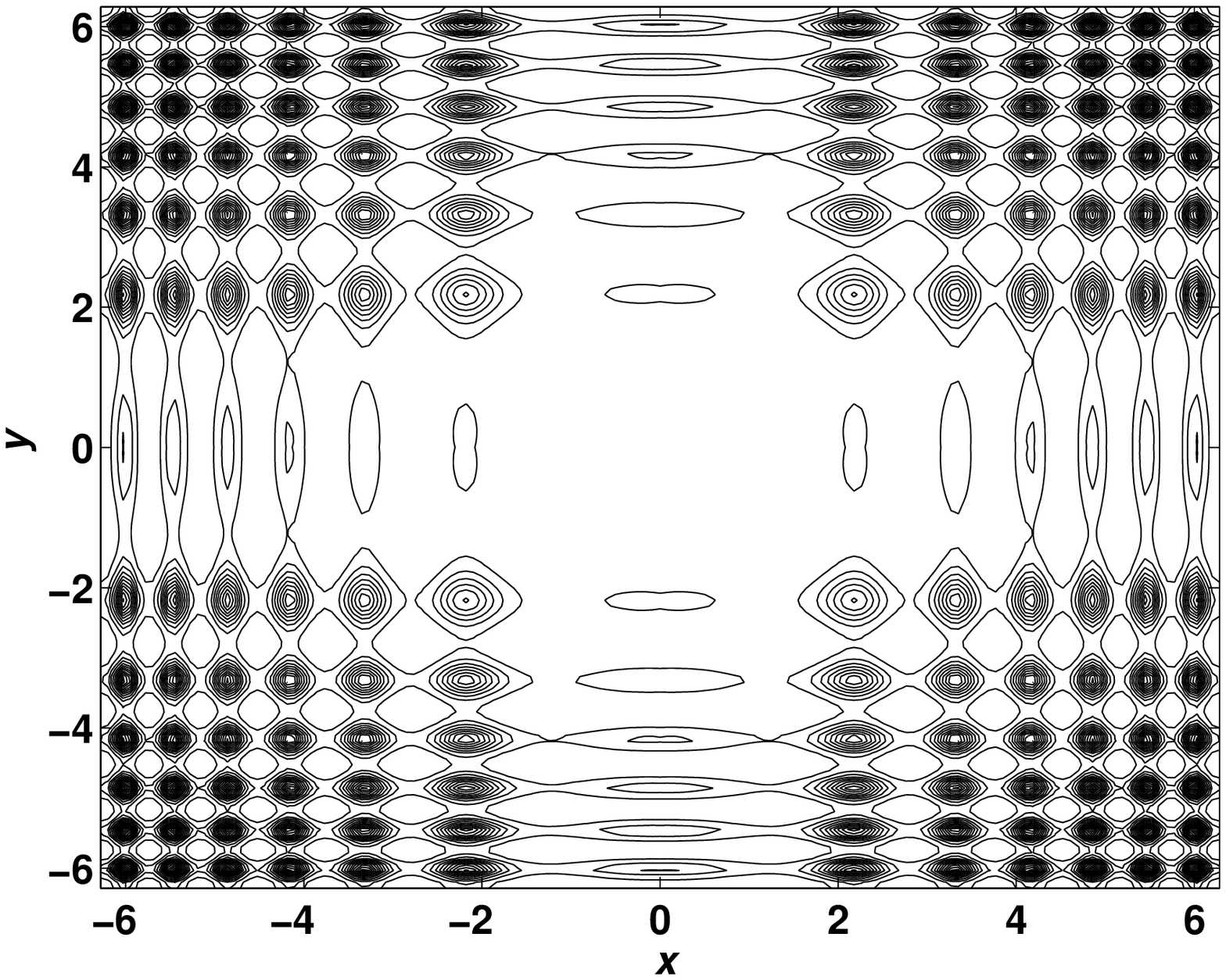} }
\caption{Yang's forest function for two independent variables
with a global minimum $f_*=0$ at $(0,0)$. \label{Yang-fig3} }
\end{figure}

\begin{figure}
 \centerline{\includegraphics[width=4in,height=3in]{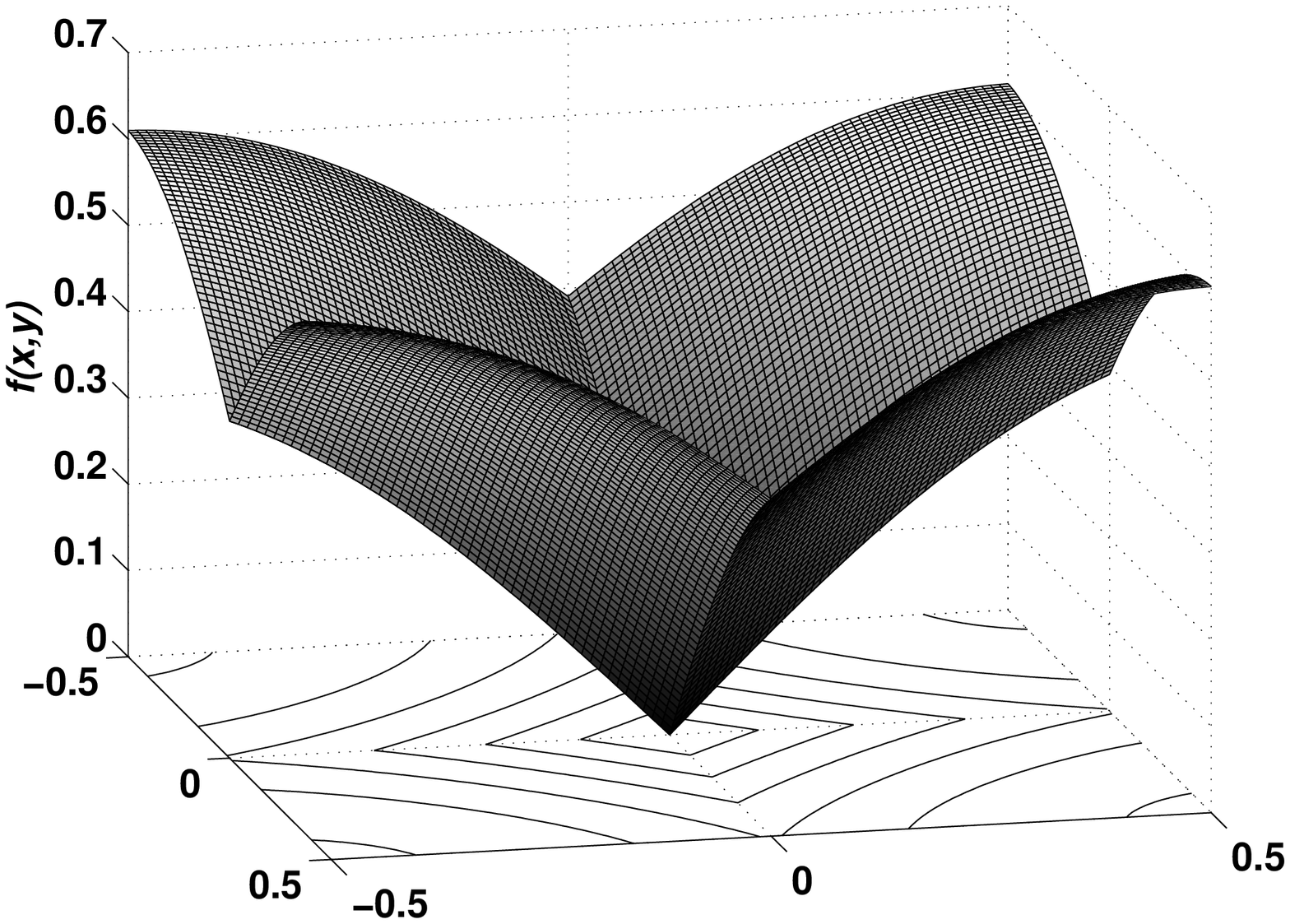}}
\caption{Non-smoothness of Yang's forest function near the global minimum $(0,0)$. \label{Yang-fig4} }
\end{figure}

\begin{table}[ht]
\caption{Comparison of algorithm performance}
\centering
\begin{tabular}{ccccc}
\hline \hline
Functions/Algorithms & GA & PSO & LFA \\
\hline
Michalewicz ($d\!\!=\!\!16$)  & $89325 \pm 7914 (95 \%)$  & $6922 \pm 537 (98\%)$  & $2889 \pm 719 (100\%)$ \\

 Rosenbrock ($d\!\!=\!\!16$) & $55723 \pm 8901 (90\%)$ & $32756 \pm 5325 (98\%)$ & $6040 \pm 535 (100\%) $ \\

 De Jong ($d\!\!=\!\!256$) & $25412 \pm 1237 (100\%)$ & $17040 \pm 1123 (100\%)$ & $5657 \pm 730 (100\%)$\\
 Schwefel ($d\!\!=\!\!128$) & $227329 \pm 7572 (95\%)$ & $14522 \pm 1275 (97\%)$ & $7923 \pm 524 (100\%)$ \\

 Ackley ($d\!\!=\!\!128$) & $32720 \pm 3327 (90\%)$ & $23407 \pm 4325 (92\%)$ & $4392 \pm 2710 (100\%)$ \\

 Rastrigin & $110523 \pm 5199 (77 \%)$ & $79491 \pm 3715 (90\%)$ & $12075 \pm 3750 (100\%)$ \\

 Easom & $19239 \pm 3307 (92\%)$ & $17273 \pm 2929 (90\%)$ & $6082 \pm 1690 (100\%)$ \\

 Griewank & $70925 \pm 7652 (90\%)$ & $55970 \pm 4223 (92\%)$ & $10790 \pm 2977 (100\%)$ \\

 Yang & $37079 \pm 8920 (88\%)$ & $19725 \pm 3204 (98\%)$ & $5152 \pm 2493 (100\%)$ \\

 Shubert (18 minima) & $54077 \pm 4997 (89\%)$ & $23992 \pm 3755 (92\%)$ & $9925 \pm 2504 (100\%)$ \\

\hline
\end{tabular}
\end{table}

We can see that the LFA is much more efficient in finding the global optima
with higher success rates. Each function evaluation is virtually instantaneous
on modern personal computer. For example, the computing time for 10,000 evaluations
on a 3GHz desktop is about 5 seconds. Even with graphics for displaying the
locations of the particles and fireflies, it usually takes less than a few minutes.
Furthermore, we have used various values of the population size $n$ or the number of
fireflies. We found that for most problems $n=15$ to $50$ would be sufficient.
For tougher problems, larger $n$ can be used, though excessively
large $n$ should not be used unless there is no better alternative,
as it is more computationally extensive.

\section{Conclusions}

In this paper, we have formulated a new L\'evy-flight firefly algorithm and analysed
its similarities and differences with particle swarm optimization. We then implemented
and compared these algorithms. Our simulation results for finding the global optima
of various test functions suggest that particle swarm often outperforms
traditional algorithms such as genetic algorithms, while LFA is
superior to both PSO and GA in terms of both efficiency and success rate.
This implies that LFA is potentially more powerful in solving NP-hard problems
which will be investigated further in future studies.

The basic L\'evy-flight firefly algorithm is very efficient. A further
improvement on the convergence of the algorithm is to carry out
sensitivity studies by varying various parameters such as $\beta_0$, $\gamma$,
$\a$ and more interestingly $\lam$.  These could form important topics for further research.
In addition, further studies on the application of FLA in combination
with other algorithms may form an exciting area for further research in optimization.

%% End of text %%

\end{document}